\documentclass[12pt]{article} 

\usepackage{epsfig,latexsym,amsfonts,amsmath,amssymb, color}
\usepackage[T1,T2A]{fontenc}
\newtheorem{theorem}{\textsf{Theorem}}[section]
\newtheorem{lemma}{\textsf{Lemma}}[section]

\newtheorem{rem}{\textsf{Remark}}
\newtheorem{defin}{\textsf{Definition}}
\newtheorem{cor}{Corollary}[section]

\def\XXint#1#2#3{{\setbox0=\hbox{$#1{#2#3}{\int}$}
	\vcenter{\hbox{$#2#3$}}\kern-.5\wd0}}

\title{A comparison theorem for nonsmooth nonlinear operators}
\author
{Vladimir Kozlov\footnote{Department of Mathematics, University of Link\"oping, SE-581 83
Link\"oping, Sweden} \ and Alexander Nazarov\footnote{St.-Petersburg Department of Steklov
Mathematical Institute, Fontanka, 27, St.-Petersburg, 191023, Russia, and
St.-Petersburg State University, Universitetskii pr. 28, St.-Petersburg, 198504, Russia}
}
\date{} 

\begin{document}
\maketitle

\begin{abstract}
We prove a comparison theorem for super- and sub-solutions with non-vanishing gradients to semilinear PDEs provided a nonlinearity $f$ is $L^p$ function with $p > 1$. The proof 
is based on a strong maximum principle for solutions of divergence type elliptic equations with VMO leading coefficients and with lower order coefficients from a Kato class. 
An application to estimation of periodic water waves profiles is given.
\end{abstract}

\section{Introduction}

Let $\Omega$ be a domain in $\mathbb R^n$, $n\geq 2$. We will consider super- and sub-solutions of the equation
\begin{equation}\label{I.1}
\Delta u+f(u)=0 \;\;\;\mbox{in $\Omega$,}
\end{equation}
where  $f$ is a real valued function from $L^p_{\rm loc}(\mathbb R) $ with some $p>1$. To make the term $f(u)$ well-defined (measurable and belonging to $L^p_{\rm loc}(\Omega)$)
we will assume that $u\in C^1(\Omega)$ and $\nabla u\neq 0$ in $\Omega$. Usually,
 f is supposed to be continuous in almost all papers,
dealing with equation (\ref{I.1}) and its generalisations (see, for
example,  \cite{Kel} and \cite{Osser} and numerous papers citing these notes).

It was shown in \cite[Remark 2.3]{Nir} that the strong maximum principal may fail if the function $f$ is only H\"older continuous with an exponent less than $1$. Optimal conditions
on smoothness of $f$ for validity of the strong maximum principal can be found in \cite{PS}. The main difference in our approach is that we compare functions in a neughborhood of a
point where the gradients are not vanishing. This allows to remove any smoothness assumptions on $f$.

One of the main results of this paper is the following assertion:

\begin{theorem}\label{T1} Let $f\in L^p_{\rm loc}(\Omega)$, $p>1$. Also let $u_1,\,u_2\in C^1_{\rm loc}(\Omega)$ have non-vanishing gradients in $\Omega$ and satisfy
the inequalities
\begin{equation}\label{I.2}
\Delta u_1+f(u_1)\geq 0 \;\;\;\mbox{and}\;\;\;\Delta u_2+f(u_2)\leq 0
\end{equation}
in the weak sense. If $u_1\leq u_2$ and $u_1(x_0)=u_2(x_0)$ for some $x_0\in\Omega$ then $u_1=u_2$ in the whole $\Omega$.
\end{theorem}
We note that the theorem is not true without assumptions that the gradients do not vanish, which follows from \cite{Nir} (see \cite{KK}).

In the case $p>n$ this theorem was proved in \cite{KK}. The proof there was based on a weak Harnack inequality for non-negative solutions to the second order equation in divergence form
\begin{equation}\label{I.3}
{\cal L}u:=D_j(a_{ji}(x)D_iu)+b_j(x)D_ju=0 \quad \mbox{in} \quad\Omega
\end{equation}
and closely connected with $L^p$ properties of the coefficients $b_j$. Therefore one of our main concerns is a strong maximum principle for solutions to (\ref{I.3}).
We always assume that the matrix $(a_{ij})$ is symmetric and uniformly elliptic:
$$
\nu|\xi|^2\le a_{ij}(x)\xi_i\xi_j\le \nu^{-1}|\xi|^2, \qquad x\in\Omega,\ \xi\in\mathbb R^n.
$$
It was proved in \cite{Tru} that if $|b|\in L^p_{\rm loc}(\Omega)$ (here $b=(b_1,\ldots,b_n)$) with $p>n$ then a non-negative weak solution to (\ref{I.3}) satisfies
(here $B_\rho(x)$ stands for the ball of radius $\rho$ centered at $x$)
\begin{equation}\label{I.3a}
\rho^{-n/\gamma}||u||_{L^\gamma(B_{2\rho}(x^0))}\leq C\inf_{x\in B_\rho(x^0)}u(x),
\end{equation}
where $B_{3\rho}(x^0)\in\Omega$ and $\gamma\in(1,n/(n-2))$. So the restriction $p>n$ in this assertion inherits in our theorem in \cite{KK}\footnote{It was pointed out in \cite[Theorem 2.5']{NU} that
(\ref{I.3a}) holds if $|b|\log^{1/2}(1+|b|)\in L^2_{\rm loc}(\Omega)$ for $n=2$ and if $|b|\in L^n_{\rm loc}(\Omega)$ for $n\ge3$.}.

For our purpose another type of assumptions on the coefficients $b_j$ are more appropriate. It is  called the Kato condition, see \cite{DH} and
\cite{Sch1}.

\begin{defin}
We say that $f\in K_{n,\alpha}$ if
\begin{equation}\label{I.4}
\sup_{x}\int\limits_{|x-y|\leq r}\frac{|f(y)|}{|x-y|^{n-\alpha}}dy\to 0 \quad\mbox{as} \quad r\to 0.
\end{equation}
\end{defin}

It was proved in \cite{Kur} that inequality (\ref{I.3a}) still holds if $|b|^2\in K_{n,2}$.
For H\"older continuous coefficients $a_{ji}$ (\ref{I.3a}) was proved in \cite{Zh} under the assumption $|b|\in K_{n,1}$.
We note that from the last assertion it follows  (\ref{I.3a}) when $|b|\in L^p_{\rm loc}$, $p>1$, depends only on one variable and
the leading coefficients are H\"older continuous.

For our applications we need the leading coefficients to be only continuous. So all above mentioned results are not enough for our purpose.
Here we prove a theorem which deals with slightly discontinuous leading coefficients and allows $b^\alpha\in K_{n,\alpha}$ with $\alpha$ close to $1$ for lower order coefficients.
In order to formulate this result we need some definitions.

\begin{defin}
Let $f(x)$  be a measurable and locally integrable function. Define a quantity
$$
f^{\#r}(x):= \frac 1{|B_r|}\int\limits_{B_r(x)}\Big| f(y)-\frac 1{|B_r|}\int\limits_{B_r(x)}\! f(z)\,dz\Big|\, dy;\qquad 
\omega_f(\rho)=\sup\limits_x\sup\limits_{r\le\rho} f^{\#r}(x).
$$
We say that $f\in VMO(\mathbb R^n)$  if $\omega_f(\rho)$ is bounded and
$\omega_f(\rho)\rightarrow0$ as $\rho\to0$. In this case the function $\omega_f(\rho)$ is called VMO-modulus of $f$.

For a bounded Lipschitz domain $\Omega$ the space $f\in VMO(\Omega)$ is introduced in the same way but the integrals in the definition of $f^{\#r}(x)$ are taken
over $B_r(x)\cap\Omega$.
\end{defin}

\begin{defin}
We say that a function $\sigma : [0,1]\rightarrow \mathbb{R}_+$ belongs to the Dini class~$\mathcal{D}$~if
$\sigma$ is increasing, $\sigma (t)/t$ is summable and decreasing.
\end{defin}
It should be noted that assumption about the decay of $\sigma(t)/t$ is not restrictive (see Remark~1.2 in \cite{AN16} for more details). We use the notation
$\varkappa(\rho)=\int\limits_0^\rho\frac {\sigma (t)}t\,dt$.

\begin{theorem}\label{T2} Let $n\geq 3$. Assume that the leading coefficients $a_{ij}\in VMO(\Omega)$. Suppose that $|b|^\alpha\in K_{n,\alpha}$  and
\begin{equation}\label{b-cond}
\sup\limits_x\int\limits_{r/2\le|x-y|\le r}\frac{|b(y)|^{\alpha}}{|x-y|^{n-\alpha}}\,dy\le \sigma^\alpha(r),
\end{equation}
for some $\alpha>1$ and $\sigma\in{\cal D}$.

If a function $u\in W^{1,p}(\Omega)$, $p>n$, satisfies $u\geq 0$ and ${\mathcal L}u\geq 0$ in $\Omega$ then either $u>0$ in $\Omega$ or $u\equiv0$ in $\Omega$.
\end{theorem}

\begin{rem}\label{alphabeta}
 Notice that the assumption $|b|^{\alpha_2}\in K_{n,\alpha_2}$ does not imply $|b|^{\alpha_1}\in K_{n,\alpha_1}$ for any $1<\alpha_1<\alpha_2$. However, if the
 condition (\ref{b-cond}) holds with $\alpha=\alpha_2$ then the H\"older inequality ensures $|b|^{\alpha_1}\in K_{n,\alpha_1}$, and (\ref{b-cond}) holds with
 $\alpha=\alpha_1$ (and another function $\sigma\in{\cal D}$).
\end{rem}

For $n=2$ we need a stronger assumption.

\begin{theorem}\label{T2X} Let $n=2$. Assume that the leading coefficients $a_{ij}\in VMO(\Omega)$. Suppose that
\begin{equation}\label{b-cond-n=2}
\sup\limits_x\int\limits_{r/2\le|x-y|\le r}\frac{|b(y)|^{\alpha}\log^{\alpha}\frac r{|x-y|}}{|x-y|^{2-\alpha}}\,dy\le \sigma^\alpha(r),
\end{equation}
for some $\alpha>1$ and $\sigma\in{\cal D}$.

If a function $u\in W^{1,p}(\Omega)$, $p>2$, satisfies $u\geq 0$ and ${\mathcal L}u\geq 0$ in $\Omega$ then either $u>0$ in $\Omega$ or $u\equiv0$ in $\Omega$.
\end{theorem}

For $\gamma\in(0,1)$ we define the annulus
\begin{equation}\label{A13s}
{\mathcal X}_{r,\gamma}(x)=\{y\in B_r(x)\,:\,|x-y|>\gamma r\}.
\end{equation}
 If the location of the center is not important we write simply $B_r$ and ${\mathcal X}_{r,\gamma}$.

 As usual, for a bounded bomain $\Omega$ we denote by $W^{1,q}_0(\Omega)$, $q>1$, the closure in $W^{1,q}(\Omega)$ of the set of smooth compactly supported function,
with the norm
$$
||u||_{W^{1,q}_0(\Omega)}=\Big(\int\limits_{\Omega}|\nabla u|^qdx\Big)^{1/q}.
$$

\section{Strong maximum principle for operators with lower order terms}

\subsection{Coercivity}

Let $\Omega'$ be a bounded subdomain in $\Omega$. Consider the problem
\begin{equation}\label{I.5}
{\mathcal L}_0u:=D_j(a_{ij}(x)D_iu)=f \quad\mbox{in}\quad \Omega';\qquad u=0\quad\mbox{on}\quad \partial \Omega'.
\end{equation}
We say that the operator ${\cal L}_0$ is $q$-coercive in $\Omega'$ for some $q>1$, if for each $f\in W^{-1,q}(\Omega')$ the problem (\ref{I.5}) has a unique solution
$u\in W^{1,q}_0(\Omega')$ and this solution satisfies
\begin{equation}\label{I.6}
||u||_{W^{1,q}_0(\Omega')}\leq C_q||f||_{W^{-1,q}(\Omega')}
\end{equation}
with $C_q$ independent on $f$ and $u$.

It is well known that for arbitrary measurable and uniformly elliptic coefficients the operator ${\cal L}_0$ is $2$-coercive in arbitrary bounded $\Omega'$. Further,
if the coefficients $a_{ij}\in VMO(\Omega)$ then the operator ${\cal L}_0$ is $q$-coercive for arbitrary $q>1$ in arbitrary bounded $\Omega'\subset\Omega$
with $\partial\Omega'\in{\cal C}^1$, see \cite{BW}. The coercivity constant $C_q$ depends on $\Omega'$ and VMO-moduli of $a_{ij}$.
Moreover, by dilation we can see that for $\Omega'=B_r$, $r\le1$, this constant does not depend on $r$. For $\Omega'={\mathcal X}_{r,\gamma}$, $r\le1$,
$C_q$ depends on $\gamma$ but not on $r$.

Let now the operator in (\ref{I.5}) be $q$-coercive for certain $q>2$. Put
$$
f=\nabla\cdot {\bf f}+f_0,\;\;
$$
where ${\bf f}=(f_1,\ldots,f_n)\in (L^q(\Omega'))^n$ and $f_0\in L^{nq/(n+q)}(\Omega')$. Then by the imbedding theorem $f\in W^{-1,q}(\Omega')$ and (\ref{I.6})
takes the form
$$
||u||_{W^{1,q}_0(\Omega')}\leq C(\sum_{j=1}^n||f_j||_{L^q(\Omega')}+||f_0||_{L^{nq/(n+q)}(\Omega')}).
$$
We need the following local estimate.

\begin{theorem}\label{coerc}
Let $\Omega'$ be a bounded Lipschitz subdomain of $\Omega$ and let the operator in (\ref{I.5}) be $q$-coercive for certain  $q>2$.
Let also $u\in W^{1,q}(\Omega')$ is such that
$$
D_j(a_{ji}(x)D_iu)=0 \quad\mbox{in} \quad\Omega'\cap B_r; \qquad u=0 \quad\mbox{on} \quad\partial \Omega'\cap B_r.
$$
Then for a fixed $\lambda\in (0,1)$
\begin{equation}\label{M12a}
||\nabla u||_{L^q(B_{\lambda r}\cap \Omega')}\leq Cr^{-1} ||u||_{L^q(B_r\cap \Omega')},
\end{equation}
where $C$ may depend on the domain $\Omega'$, $q$, $\lambda$ and the coercivity constant $C_q$ but it is independent of $r$.
\end{theorem}
{\bf Proof.} First, we claim that the problem (\ref{I.5}) is $s$-coercive for any $s\in [2,q]$. Indeed, we have coercivity for $s=2$ and $s=q$, and the claim follows
by interpolation.

Second, the estimate (\ref{M12a}) for $q=2$ follows by the the Caccioppoli inequality. For $q>2$ we choose a cut-off function $\zeta$
such that $\zeta=1$ on $B_{\lambda r}$ and $\zeta=0$ outside $B_{\lambda_1r}$, where $\lambda<\lambda_1<1$, and $\nabla\zeta \leq cr^{-1}$. Then
$$
D_j(a_{ji}(x)D_i(\zeta u))=D_j(a_{ji}(x)uD_i\zeta )+(D_j\zeta)a_{ji}(x)D_i u \quad\mbox{in} \quad\Omega'; \qquad \zeta u=0 \quad\mbox{on} \quad\partial \Omega'.
$$
Then by the $s$-coercivity of the operator we have
$$
||\nabla u||_{L^s(B_{\lambda r})}\leq Cr^{-1}(||u||_{L^s(B_{\lambda_1}r)}+||\nabla u||_{L^{ns/(n+s)}(B_{\lambda_1r})})
$$
We choose $s=\min (q,2n/(n-2))$. Then the last term in the right-hand side is estimated by $Cr^{n/s-n/2}||\nabla u||_{L^{2}(B_{\lambda_1r})}$,
and hence by the proved estimate for $q=2$, we obtain (\ref{M12a}) for $q=s$. Repeating this argument (but using now the estimate (\ref{M12a}) for $q=s$)
we arrive finally at (\ref{M12a}).

\subsection{Estimates of the Green functions}

Let $\cal L$ be an operator of the form (\ref{I.3}), and the assumptions of Theorem \ref{T2} are fulfilled.
We establish the existence and some estimates of the Green function $G=G(x,y)$ for the problem
\begin{equation}\label{A11a}
{\mathcal L}u=f\quad\mbox{in}\quad B_r;\qquad u=0\quad\mbox{on}\quad \partial B_r,
\end{equation}
in sufficiently small ball $B_r\subset \Omega$.

\begin{lemma}\label{LA1a} Let $n\geq 3$.
There exists a positive constant $R$ depending on $n$, the ellipticity constant $\nu$, VMO-moduli of coefficients $a_{ij}$,
the exponent $\alpha$ and the function $\sigma$ from (\ref{b-cond}) such that for any $r\le R$ there is the Green function $G(x,y)$ of the problem (\ref{A11a})
in a ball $B_r\subset\Omega$. Moreover, it is continuous w.r.t. $x$ for $x\ne y$ and satisfies the estimates
\begin{eqnarray}\label{A12a}
&&|G(x,y)|\le\frac{C_1}{|x-y|^{n-2}}; \nonumber\\
&\mbox{if}&|x-y|\le {\rm dist}(x,\partial B_r)/2 \quad  \quad\mbox{then}\quad
G(x,y)\ge\frac{C_2}{|x-y|^{n-2}},
\end{eqnarray}
where the constants $C_1$ and $C_2$ depend on the same quantities as $R$.
\end{lemma}

{\bf Proof}. We use the idea from \cite{Zh}.
Denote by $G_0(x,y)$ the Green function of the problem (\ref{I.5}) in the ball $B_r$.
The estimates (\ref{A12a}) for $G_0$ were proved in \cite{LS} (see also \cite{GW}): 
\begin{eqnarray}\label{A13a}
&0<&G_0(x,y)\le\frac{C_{10}}{|x-y|^{n-2}}; \nonumber\\
&\mbox{if}& |x-y|\le {\rm dist}(x,\partial B_r)/2 \quad  \quad\mbox{then}\quad
G_0(x,y)\ge\frac{C_{20}}{|x-y|^{n-2}},
\end{eqnarray}
where the constants $C_{10}$ and $C_{20}$ depend only on $n$ and $\nu$.

By Remark \ref{alphabeta}, we can assume without loss of generality $\alpha< n/(n-1)$. Put $q=\alpha'> n$ and denote by $C_q$ the coercivity
constant for ${\cal L}_0$ in the ball. We begin with the estimate for any $B_\rho(y)\subset\Omega$
\begin{eqnarray*}
&&\int\limits_{B_\rho(y)}\,|b(x)|\, |D_xG_0(x,y)|\,dx
=\sum\limits_{j=0}^{\infty}
\,\int\limits_{{\cal X}_{\frac {\rho}{2^j},\frac 12}(y)}|b(x)|\, |D_xG_0(x,y)|\,dx\\
&\le&\sum\limits_{j=0}^{\infty} \Big(\int\limits_{{\cal X}_{\frac {\rho}{2^j},\frac 12}(y)}|D_xG_0(x,y))|^q\,dx\Big)^{\frac 1q}
\Big(\int\limits_{{\cal X}_{\frac {\rho}{2^j},\frac 12}(y)}|b(x)|^{q'}\,dx\Big)^{\frac 1{q'}}=:\sum\limits_{j=0}^{\infty} A_{j1}A_{j2}.
\end{eqnarray*}
By (\ref{M12a}) and (\ref{A13a}), we have $A_{j1}\le C(n,\nu,C_q)(2^{-j-1}\rho)^{1- \frac n{q'}}$, and (\ref{b-cond}) gives
$$
A_{j1}A_{j2}\le 2^{\frac n{q'}-1}C\,\Big(\int\limits_{{\cal X}_{\frac {\rho}{2^j},\frac 12}(y)}\frac {|b(x)|^{q'}}{|x-y|^{n- q'}}\,dx\Big)^{\frac 1{q'}}
\le 2^{\frac n{q'}-1}C\sigma(2^{-j}\rho).
$$
Therefore,
\begin{equation}\label{varkappa}
\int\limits_{B_\rho(y)}\,|b(x)|\, |D_xG_0(x,y)|\,dx\le 2^{\frac n{q'}-1}C\sum\limits_{j=0}^{\infty}\sigma(2^{-j}\rho)
\le A(n,\nu,q,C_q)\int\limits_0^\rho \frac {\sigma (t)}t\,dt\equiv A\varkappa(\rho).
\end{equation}

Next, we write down the equation for $G$
\begin{equation}\label{series}
{\cal L}_0(G-G_0)=-b\cdot DG \qquad \Longleftrightarrow \qquad (I+G_0*(b\cdot D))(G-G_0)=-b\cdot DG_0
\end{equation}
and obtain
$$
G=G_0-G_0*(b\cdot DG_0)+G_0*(b\cdot DG_0)*(b\cdot DG_0)-\dots =:\sum_{k=0}^{\infty}J_k
$$
provided this series converges.

We claim that
\begin{equation}\label{J-est}
|J_k(x,y)|\le\frac{C_{10}C^k\varkappa^k(r)}{|x-y|^{n-2}}\quad \Longrightarrow \quad
|J_{k+1}(x,y)|\le\frac{C_{10}C^{k+1}\varkappa^{k+1}(r)}{|x-y|^{n-2}},
\end{equation}
for a proper constant $C$. Indeed,
$$
J_{k+1}(x,y)=\int\limits_{B_r}J_k(x,z)(b(z)\cdot D_zG_0(z,y))\,dz.
$$
Denote $2\rho=|x-y|$. Then
$$
|J_{k+1}(x,y)|\le\Big(\int\limits_{B_r\setminus B_\rho(x)}+\int\limits_{B_r\cap B_\rho(x)}\Big)|J_k(x,z)|\,|b(z)|\, |D_zG_0(z,y)|\,dz=: I_1+I_2.
$$
We have by (\ref{varkappa})
\begin{eqnarray*}
I_1 &\le& \int\limits_{B_r\setminus B_\rho(x)}
\frac{C_{10}C^k\varkappa^k(r)}{|x-z|^{n-2}}\,|(b(z)|\, |D_zG_0(z,y))|\,dz\\
&\le& \frac{2^{n-2}C_{10}C^k\varkappa^k(r)}{|x-y|^{n-2}}\int\limits_{B_r\setminus B_\rho(x)}|b(z)|\, |D_zG_0(z,y)|\,dz
\le\frac{2^{n-1}AC_{10}C^k\varkappa^{k+1}(r)}{|x-y|^{n-2}}
\end{eqnarray*}
(here we used an evident inequality $\varkappa(2r)\le 2\varkappa(r)$).
Further,
$$
I_2\le \int\limits_{B_r\cap B_\rho(x)} \frac{AC^k\varkappa^k(r)}{|x-z|^{n-2}}\,|b(z)|\,|D_zG_0(z,y)|\,dz.
$$
By (\ref{M12a}) and (\ref{A13a}),
$$
\Big(\int\limits_{B_r\cap B_\rho(x)}|DG_0(z,y))|^q\,dy\Big)^{\frac 1q}\le A\rho^{1-\frac n{q'}}.
$$
Using the assumption $q>n$ we get
$$
I_1\le \frac {A^2C^k\varkappa^k(r)}{\rho^{\frac n{q'}-1}}\,
\Big(\int\limits_{B_r\cap B_\rho(x)}\frac{|b(z)|^{q'}}{|x-z|^{(n-2)q'}}\,dz\Big)
^{\frac 1{q'}}\le \frac {A^2C^k\varkappa^k(r)}{\rho^{\frac n{q'}-1}}\,
\rho^{1-\frac nq}\varkappa(\rho)
\le\frac{2^{n-2}A^2C^k\varkappa^{k+1}(r)}{|x-y|^{n-2}},
$$
and the claim follows if we put $C=2^{n-2}(A+2C_{10})$.

Thus, the series in (\ref{series}) converges if $\varkappa(r)$ is sufficiently small. Moreover, if $2^{n-2}(A+2C_{10})\varkappa(r)\le \frac {C_{20}}{C_{20}+2C_{10}}$
then (\ref{A13a}) implies (\ref{A12a}) with $C_1=C_{10}+\frac {C_{20}}2$, $C_2=\frac {C_{20}}2$. \medskip

To prove the continuity of $G$ we take $\widehat x$ such that $|x-\widehat x|\le \rho/2=|x-y|/4$. Since $q>n$, the estimate (\ref{M12a}) and the Morrey embedding theorem give
$$
|G_0(x,y)-G_0(\widehat x,y)|\le C_{30}(n,\nu,C_q)\,\frac {|x-\widehat x|^{1-\frac nq}}{|x-y|^{\frac n{q'}-1}}.
$$
We write down the relation
$$
G(x,y)-G_(\widehat x,y)=\sum_{k=0}^{\infty}\Big(J_k(x,y)-J_k(\widehat x,y)\Big)
$$
and deduce, similarly to (\ref{J-est}), that
\begin{multline*}\label{J-J-est}
|J_k(x,y)-J_k(\widehat x,y)|\le C_{30}C^k\varkappa^k(r)\,\frac{|x-\widehat x|^{1-\frac nq}}{|x-y|^{\frac n{q'}-1}}\quad \Longrightarrow\\
\Longrightarrow \quad
|J_{k+1}(x,y)-J_{k+1}(\widehat x,y)|\le C_{30}C^{k+1}\varkappa^{k+1}(r)\,\frac{|x-\widehat x|^{1-\frac nq}}{|x-y|^{\frac n{q'}-1}}.
\end{multline*}
Therefore, if $\varkappa(r)$ is sufficiently small,
$$
|G(x,y)-G(\widehat x,y)|\le C_3\,\frac {|x-\widehat x|^{1-\frac nq}}{|x-y|^{\frac n{q'}-1}}.
$$
\hfill$\square$\medskip

\begin{rem}
 In fact, since $q$ can be chosen arbitrarily large, $G$ is locally H\"older continuous w.r.t. $x$ for $x\ne y$ with arbitrary exponent $\beta\in(0,1)$.
\end{rem}

Now let ${\cal X}_{r,\gamma}\subset\Omega$. Consider the Dirichlet problem
\begin{equation}\label{A12c}
{\mathcal L}u=f\quad\mbox{in}\quad {\cal X}_{r,\gamma};\qquad u=0\quad\mbox{on}\quad \partial {\cal X}_{r,\gamma}.
\end{equation}

\begin{lemma}\label{LA12d}
The statement of Lemma \ref{LA1a} holds for the problem (\ref{A12c}). The constants $R$, $C_1$ and $C_2$ may depend on the same quantities as in Lemma \ref{LA1a}
and also on $\gamma$.
\end{lemma}

The proof of Lemma \ref{LA1a} runs without changes.

\subsection{Approximation Lemma and weak maximum principle}

\begin{lemma}\label{LA12e} Under assumptions of Lemma \ref{LA1a}, let $u$ be a weak solution of the equation
$$
{\cal L}u=f \quad {\rm in}\quad B_r\subset\Omega, 
$$
with $r\le R$ where $R$ is the constant from Lemma \ref{LA1a}. Let $f$ be a finite signed measure in $B_r$.

Put $b^m_i:=b_i\chi_{\{|b|\leq m\}}$ 
and define a sequence $f_m\in L^{\infty}(B_r)$ such that $f_m\to f$ in the space of measures.
Denote by $u_m
$ the solution of the problem
$$
{\cal L}_mu_m:=-D_i(a_{ij}D_ju_m)+ b^m_iD_iu_m=f^m\quad\mbox{in}\quad B_r;\qquad u^m=u\quad\mbox{on}\quad \partial B_r.
$$
Then
$$
\int\limits_{B_r}|u_m(x)-u(x)|\,dx\rightarrow 0\quad\mbox{as}\quad m\rightarrow\infty.
$$	
\end{lemma}

{\bf Proof.} It is easy to see that the difference $v_m=u_m-u$ solves the problem
$$
{\cal L}_m v_m=
(b_i-b_i^m)D_iu+f_m-f\quad\mbox{in}\quad B_r; \qquad v_m=0\quad\mbox{on}\quad \partial B_r.
$$
Using the Green function $G_m$ of the operator ${\cal L}_m$ in $B_r$ with the Dirichlet boundary conditions we
get
\begin{eqnarray*}
\int\limits_{B_r}|v_m(x)|\,dx\le
\Big(\int\limits_{B_r}|(b_i-b_i^m)D_iu| + \int\limits_{B_r}|f_m-f|\Big)
\cdot\sup\limits_y\int\limits_{B_r}|G_m(x,y)|\,dx.
\end{eqnarray*}
By Lemma \ref{LA1a}, $|G_m(x,y)|\le C|x-y|^{2-n}$ 
with constant independent of $m$. Thus, the supremum of the last integral is bounded. The first integral in brackets tends to zero by the Lebesgue Dominated convergence Theorem, and
the Lemma follows.\hfill$\square$\medskip

\begin{cor}\label{K1a}
If ${\cal L}u=f
\geq 0$ in $B_r$,
and $u\geq 0$ on $\partial B_r$, then $u\geq 0$ in $B_r$.
\end{cor}

This statement follows from standard weak maximum principle and Lemma \ref{LA12e}.

\begin{lemma}\label{LA12f}
Let ${\cal X}_{r,\gamma}\subset\Omega$ and let $r\le R$ where $R$ is the constant from Lemma \ref{LA12d}. Then the assertion of
Lemma \ref{LA12e} and the corollary \ref{K1a} are still true.
\end{lemma}

\subsection{Strong maximum principle}


\begin{lemma}\label{LA13a} Let $B_r\subset \Omega$ be a ball and let $r\le R$ where $R$ satisfies the assumptions of  Lemma \ref{LA1a} and Lemma \ref{LA12d} with $\gamma=1/4$.
Then the Green function of $\cal L$ in $B_r$ is strictly positive: $G(x,y)>0$ for $x,y\in B_r$, $x\neq y$.
\end{lemma}

{\bf Proof.}
 First suppose that $G(x^*,y)<0$ for certain $x^*\in B_r$ and for a positive measure set of $y$. By continuity of $G$ in $x$ we have $G(x,y)<0$ for
 a (maybe smaller) positive measure set of $y$ and an open set of $x$. Therefore, we can choose a bounded nonnegative function $f$ such that
$u(x)=\int\limits_{B_r}G(x,y)f(y)<0$ on an open set. But this would contradict to the weak maximum principle, see Corollary \ref{K1a}. Thus, we can change
$G$ on a null measure set and assume it nonnegative.

Next, suppose that for certain $y^*$ the set ${\cal S}(y^*)=\{x\in B_r\,:\,G(x,y^*)=0\}$ is non-empty.
We choose then $x^0\in {\cal S}(y^*)$ and $y^0\in B_r\setminus {\cal S}(y^*)$ such that $\rho:=|x^0-y^0|={\rm dist} (y^0, {\cal S}(y^*))$.
Due to the second estimate in (\ref{A12a}) $x^0$ is separated from $y^*$, while $\rho$ can be chosen arbitrarily small. So, we can suppose that
$$
B_{2\rho}(y^0)\subset B_r\setminus \{y^*\};\qquad B_{\rho/2}(y^0)\cap {\cal S}(y^*)=\emptyset.
$$

We introduce the Green function $\widehat{G}$ of $\cal L$ in $B_{2\rho}(y^0)$ and claim that the function $\delta\widehat{G}(\cdot, y^0)$ with sufficiently small
$\delta>0$ is a lower barrier for $G(\cdot,y^*)$ in the annulus ${\mathcal X}_{2\rho,1/4}(y^0)$. Indeed, the boundary $\partial{\mathcal X}_{2\rho,1/4}(y^0)$
consists of two spheres. Notice that $G(\cdot,y^*)>0$ on $\partial B_{\rho/2}(y^0)$ while $\widehat{G}(\cdot,y^0)=0$ on $\partial B_{2\rho}(y^0)$. Thus,
there exists a positive $\delta$ such that
$$
G(\cdot,y^*)\geq \delta \widehat{G}(\cdot,y^0)\qquad\mbox{on}\quad \partial {\mathcal X}_{2\rho,1/4},
$$
and the claim follows. By Lemma \ref{LA12f} $G(\cdot,y^*)\geq \delta \widehat{G}(\cdot,y^0)$ in the whole annulus, and, in particular,
$$
G(x^0,y^*)\geq \delta \widehat{G}(x^0,y^0)>0
$$
(the last inequality follows from the second estimate in (\ref{A12a})). The obtained contradiction proves the Lemma. \hfill$\square$\medskip

{\bf Proof of Theorem \ref{T2}}. We repeat in essential the proof of Lemma \ref{LA13a}. Denote the set ${\cal S}=\{x\in\Omega\,:\, u(x)=0\}$ and suppose that ${\cal S}\neq \Omega$.
Then we can choose $x^0\in {\cal S}$ and $y^0\in \Omega\setminus {\cal S}$ such that $\rho:=|x^0-y^0|={\rm dist} (y^0, {\cal S})$, and $\rho$ can be chosen arbitrarily small.
Repeating the proof of Lemma \ref{LA13a} we introduce the same Dirichlet Green function $\widehat{G}(\cdot, y^0)$ of $\cal L$ in $B_{2\rho}(y^0)$ and show that
$\delta\widehat{G}(\cdot, y^0)$ with sufficiently small $\delta>0$ is a lower barrier for $u$ in the annulus ${\mathcal X}_{2\rho,1/4}(y^0)$. This ends the proof. \hfill$\square$

\subsection{The case $n=2$}

The case $n=2$ is treated basically in the same way as the case $n\geq 3$, but some changes must be done mostly due to the fact that the estimate of the Green function contains logarithm.

Let us explain what changes must be done in the argument in compare with $n\geq 3$.
Denote by $G_0(x,y)$ the Green function of the problem (\ref{I.5}) in the disc $B_r$. Then for $x\ne y\in B_r$,
\begin{eqnarray}\label{Kx1}
&0<&G_0(x,y)\le C'_{10}\log \Big(\frac r{|x-y|} + 2\Big); \nonumber\\
&\mbox{if}& |x-y|\le {\rm dist}(x,\partial B_r)/2 \quad  \quad\mbox{then}\quad
G_0(x,y) \ge C'_{20}\log \Big(\frac r{|x-y|} + 2\Big),
\end{eqnarray}
where the constants $C'_{10}$ and $C'_{20}$ depend only on the ellipticity constants of the operator ${\mathcal L}_0$. Indeed by \cite{LS} these estimates can be reduced to similar estimates
for the Laplacian, when they can be verified directly (in this case the Green function can be written explicitly).

Analog of Lemma \ref{LA1a} reads as follows.

\begin{lemma}\label{LA1aX} Let $n=2$.
There exists a positive constant $R$ depending on the ellipticity constant $\nu$, VMO-moduli of coefficients $a_{ij}$,
the exponent $\alpha$ and the function $\sigma$ from (\ref{b-cond-n=2}) such that for any $r\le R$ there is the Green function $G(x,y)$ of the problem (\ref{A11a})
in a disc $B_r\subset\Omega$. Moreover, it is continuous w.r.t. $x$ for $x\ne y$ and satisfies the estimates
\begin{eqnarray}\label{A12aX}
&&|G(x,y)|\le C'_1\log \Big(\frac r{|x-y|} + 2\Big); \nonumber\\
&&\mbox{if}\quad |x-y|\le {\rm dist}(x,\partial B_r)/2 \quad  \quad\mbox{then}\quad
|G(x,y)|\ge C'_2\log \Big(\frac r{|x-y|} + 2\Big),
\end{eqnarray}
where the constants $C'_1$ and $C'_2$ depend on the same quantities as $R$.
\end{lemma}
To prove this statement we establish the inequality (\ref{varkappa}) by using the estimate (\ref{Kx1}) and the assumption (\ref{b-cond-n=2}) instead of (\ref{b-cond}). The rest of the proof
runs without changes.

Corresponding analog of Lemma \ref{LA12d} is true here also.

The remaining part of the proof of Theorem \ref{T2X} is the same as that of Theorem \ref{T2}.

\section{Comparison theorem for nonlinear operators}

\subsection{Proof of Theorem \protect\ref{T1}}

We recall that the statement of Theorem \ref{T1} with $p>n$ was proved in \cite[Theorem 1]{KK} by reducing it to the strong maximum principle for the equation (\ref{I.3}) with continuous
leading coefficients and $f(x_1)$ playing the role of a coefficient $b_1$. Since the function $f$ depends only on one variable, the assumption $f\in L^p_{\rm loc}(\mathbb R)$ with a certain $p>1$
implies (\ref{b-cond}) and (\ref{b-cond-n=2}) with $\alpha=p$. Thus, we can apply our Theorems \ref{T2} and \ref{T2X} (for $n\ge3$ and $n=2$, respectively) instead of \cite[Lemma 1]{KK},
and the proof of \cite[Theorem 1]{KK} runs without other changes.

\subsection{Application to water wave theory}

In the paper \cite{KK}, two theorems were proved on estimates of the free surface profile of water waves on two-dimensional flows with vorticity in a channel, see
Theorems 2 and 3 in \cite{KK}. The vorticity function $\omega$ in that theorems was assumed to belong to $L^p_{\rm loc}(\mathbb R)$ with $p>2$ and the proof was based on the application of
\cite[Theorem 1]{KK}. Now the application of our Theorem \ref{T1} 
allows us to weaken the apriori assumption for the vorticity function to $\omega\in L^p_{\rm loc}(\mathbb R)$, $p>1$.

\vspace{6mm}

\noindent {\bf Acknowledgements.}  V.~K. acknowledges the support of the Swedish Research Council (VR) Grant EO418401. A.~N. was partially supported by Russian Foundation for Basic Research,
Grant 18-01-00472. He also thanks the Link\"oping University for the hospitality during his visit in January 2018.

\end{document}